
\magnification=\magstep1
\input amstex
\documentstyle{amsppt}
\TagsOnRight
\NoRunningHeads
\newdimen\notespace \notespace=7mm
\newdimen\maxnote  \maxnote=13mm
\define\marg{\strut\vadjust{\kern-\dp\strutbox
  \vtop to\dp\strutbox{\vss \baselineskip=\dp\strutbox
  \moveleft\notespace\llap{\hbox to \maxnote{\hfil$\to$}}\null}}\ }

\redefine\marg{\relax\ }

\define\obs{\strut\vadjust{\kern-\dp\strutbox
  \vtop to\dp\strutbox{\vss \baselineskip=\dp\strutbox
  \moveleft\notespace\llap{\hbox to \maxnote{\hfil$!$}}\null}}\ }

\define\local{k[x]_{(x)}}
\define\Spec{\operatorname{Spec}}
\define\Hom{\operatorname{Hom}}
\topmatter
\title The Hilbert scheme\endtitle
\title The Hilbert scheme\\
 parameterizing finite length subschemes of the line\\
 with support at the origin
\endtitle 
\author Dan Laksov and Roy M. Skjelnes\endauthor
\affil Department of Mathematics, KTH\endaffil
\address KTH, S-100 44 STOCKHOLM, Sweden\endaddress
\email laksov\@math.kth.se, skjelnes\@math.kth.se\endemail

\abstract We introduce symmetrizing operators of the polynomial ring
$A[x]$ in the variable $x$ over a ring $A$. When $A$ is an algebra
over a field $k$ these operators are used to characterize the monic
polynomials $F(x)$ of degree $n$ in $A[x]$ such that $A\otimes_k
\local/(F(x))$ is a free $A$--module of rank $n$. We use the
characterization to determine the Hilbert scheme parameterizing
subschemes of length $n$ of $\local$.\endabstract
\endtopmatter
\document

\subhead Introduction\endsubhead We shall study the Hilbert scheme
parameterizing finite length subschemes of the local ring $\local$ of
the line at the origin. The Hilbert schemes parameterizing finite
length subschemes of local rings have mostly been studied for local
rings at smooth points on surfaces (see e.g. \cite{B}, \cite{BI},
\cite{C}, \cite{G}, \cite{I1}, \cite{I2}, \cite{I3}, \cite{P}). The
focus has been on the rational points of the Hilbert schemes rather
than on the schemes themselves.

The purpose of the following work is to point out that we loose
essential information about the Hilbert schemes parameterizing finite
length subschemes of a local ring by considering rational points
instead of families. Indeed, there is only one rational point
$k[x]/(x^n)$ of the Hilbert scheme parameterizing subschemes of
$\local$ of length $n$ whereas, as we shall show in this article, the
Hilbert scheme is affine of dimension $n$. The coordinate ring is equal
to the localization of the symmetric polynomials of the ring $k[t_1,
\dots , t_n]$ in $n$ variables, in the multiplicatively closed subset
consisting of the products $g(t_1)\cdots g(t_n)$ for all polynomials
$g(x)$ in one variable over $k$ such that $g(0) \neq 0$. 

In forthcoming work \cite{S1} the second author will use the
techniques and results of the present article to show that the {\it
functor of families with support at the origin}, in contrast to the
Hilbert functor, is not even representable. The functor of families
with support at the origin is frequently used by some authors because
it has the same rational points as the Hilbert scheme. In \cite{S2}
the second author shows how the techniques of the present article can
be used on any localization of $k[x]$, over any ring $k$, and gives the relation to the
well known result that the Hilbert scheme of the projective line is
given by the symmetric product.

An easy and fundamental result in commutative algebra states that the
residue ring $A[x]/(F(x))$ of the polynomial ring in a variable $x$
over a ring $A$ by the ideal generated by a monic polynomial $F(x)$ of
degree $n$, is a free $A$--module of rank $n$. The key to the study of
the Hilbert scheme parameterizing finite length subschemes of $\local$
is to determine when the ring $A\otimes_k \local /(F(x))$ is a free
$A$--module of rank $n$. Easy examples show that the $A$--module
$A\otimes_k \local/(F(x))$ neither has to be of rank $n$ nor to be
finitely generated. Indeed, we have that $A\otimes_k \local/(x-1)=0$
for all $A$, and $k[u]\otimes_k \local/(x-u)=k[u]_{(u)}$ when $A=k[u]$
is a polynomial ring over $k$ in the variable $u$.

Theorem (2.3),\marg which is the main result of the article,
characterizes the monic polynomials $F(x)$ in $A[x]$ of degree $n$
such that $A\otimes_k \local/(F(x))$ is a free $A$--module of rank
$n$. The essential technical tool used in the proof is the
introduction of symmetrizing operators on the polynomials in the ring
$A[x]$ associated to $F(x)$. The method is introduced in Section 1 and
is the main technical novelty of the article.

Theorem (2.3)\marg is used in Section 3, where we describe the ideals
$I$ in $A\otimes_k \local$ such that $A\otimes_k \local/I$ is a free
$A$--module of rank $n$. We then proceed in Section 4 to determine the
Hilbert scheme parameterizing length $n$ subschemes of $\local$.

\subhead 1. Notation and the symmetrizing operators\endsubhead

\definition{1.1. Notation}Given a commutative ring $A$. Denote by
$A[x,t] =A[x, t_1, \dots , t_n]$ the polynomial ring over $A$ in the
variables $x, t_1, \dots , t_n$. Let $\varphi\colon A\to K$ be a
homomorphism of commutative rings. We shall consider $K$ as an
$A$--algebra via this homomorphism and write $K[x,t] =K\otimes_A
A[x,t]$. Given a polynomial $G(x) =g_0x^m +\cdots +g_m$ in $A[x]$, we
write $G^\varphi(x) =\varphi(g_0) x^m +\cdots +\varphi(g_m)$ in
$K[x]$.  \enddefinition

We denote by $s_i(t)$ the $i$'th elementary symmetric function in
the variables $t_1, \dots , t_n$.

\definition{1.2. The main construction}Given a polynomial $G(x)$ in
$A[x]$, we write 
$$s_i(G(t)) =s_i(G(t_1), \dots , G(t_n)).$$ The polynomial $s_i(G(t))$
is symmetric in the variables $t_1, \dots , t_n$. We note that the
symmetric function $s_i(x(t))$ associated to the polynomial $G(x)=x$
is equal to the elementary symmetric function $s_i(t)$ so there is no
confusion of notation. We write
 $$\multline
\Delta (G,t) =\prod_{i=1}^n \left( G(x) -G(t_i)\right) \\
 =G(x)^n  -s_1(G(t)) G(x)^{n-1} +\cdots +(-1)^n
s_n(G(t))\endmultline\tag{1.2.1}$$
 in $A[x,t]$. The polynomial  $\Delta (G,t)$ is symmetric in the
 variables $t_1,  \dots , t_n$ and $\Delta(x,t) =\prod_{i=1}^n (x-t_i)
 =x^n -s_1(t)  x^{n-1} +\cdots +(-1)^n s_n(t)$. Since $G(x) -G(t_i)$
is divisible by $x-t_i$ we obtain
 that
 $$\Delta (G,t) =H(x,t) \Delta (x,t)\tag{1.2.2}$$
in $A[x, t]$.

Fix a polynomial 
 $$F(x) =x^n -u_1x^{n-1} +\cdots +(-1)^n u_n$$
 in $A[x]$. There is a unique $A$--algebra homomorphism
 $$u\colon A[s_1(t), \dots , s_n(t)] \to A$$
 determined by $u(s_i(t)) =u_i$ for $i=1, \dots , n$. We have that 
 $\Delta^u (x,t) =F(x)$. Write  $s_{F,i}(G(t))= u(s_i(G(t)))$.  It
 follows from the formulas  (1.2.1)\marg and (1.2.2)\marg that 
 $$G(x)^n-s_{F,1}(G(t)) G(x)^{n-1} +\cdots +(-1)^n s_{F,n}(G(t))
 = H^u(x)F(x).\tag{1.2.3}$$
 in $A[x]$.
\enddefinition

\proclaim{1.3. Lemma}Given a polynomial $G(x)$ in $A[x]$. If
$s_{F, n}(G(t))$ is invertible in $A$ we have that the class of
$G(x)$ in $A[x]/(F(x))$ is invertible.
\endproclaim

\demo{Proof}When $s_{F, n}(G(t))$ is invertible we obtain from
formula (1.2.3)\marg the formula
 $$\multline
(-1)^{n+1}s_{F,n}(G(t))^{-1} G(x)\left[ G(x)^{n-1} -s_{F,1}(G(t))
   G(x)^{n-2}\right.+\cdots\\ \left. +(-1)^{n-1} s_{F,
n-1}(G(t))\right]   =1 +(-1)^{n+1}
s_{F,n}(G(t))^{-1}H^u(x)F(x).\endmultline$$  
 The Lemma follows immediately from the latter formula.
\enddemo

\proclaim{1.4. Lemma}Given a ring homomorphism $\varphi\colon A\to K$
of the ring $A$ into a field $K$. Let $\alpha_1, \dots , \alpha_n$ be
the roots of the polynomial
 $$F^\varphi(x) =x^n -\varphi(u_1) x^{n-1} +\cdots
 +(-1)^n\varphi(u_n)$$
 in the algebraic closure of $K$. Then we have that
 $$\varphi(s_{F,n}(G(t))) =G^\varphi(\alpha_1) \cdots
 G^\varphi(\alpha_n)$$
 in ${K}$.

In particular, if $k$ is a field and $\varphi$ is a homomorphism of
$k$--algebras, we have for each polynomial $g(x)$ in $k[x]$ that
 $$\varphi\left(s_{F,n}(g(t))\right) =g(\alpha_1)\cdots g(\alpha_n)$$
 in ${K}$.
\endproclaim

\demo{Proof}From the construction of Section (1.2)\marg for the ring
$K$ we obtain an expression
 $$\Delta(G^\varphi, t) =G^\varphi(x)^n -s_1(G^\varphi(t))
 G^\varphi(x)^{n-1} +\cdots +(-1)^n s_n(G^\varphi(t))$$
 in $K[x,t]$. Moreover, from the polynomial $F^\varphi(x)$ we obtain a
 unique $K$--algebra homomorphism
 $$\varphi u\colon K[s_1(t), \dots , s_n(t)] \to K$$
 determined by $(\varphi u)(s_i(t)) =\varphi(u_i)$. It follows from the
 construction of Section (1.2)\marg applied to $A$ and to $K$, that
 we have
 $$\varphi(s_{F,i}(G(t))) =s_{F^{\varphi},i} (G^\varphi(t))$$
 for $i=1, \dots ,n$. Denote by $\overline{K}$ the algebraic closure
 of $K$. The $K$--algebra homo\-mor\-phism
 $$\alpha\colon K[t_1, \dots , t_n]\to \overline{K}$$
 determined
 by $\alpha(t_i) =\alpha_i$ extends the homomorphism $\varphi u$
 because $\alpha(s_i(t)) =s_i(\alpha_1, \dots , \alpha_n)
 =\varphi(u_i) =(\varphi u)(s_i(t))$. We have that $s_n(G^\varphi(t))
=G^\varphi(t_1)\cdots 
 G^\varphi(t_n)$ in $K[t]$. Hence
 $$\varphi\left(s_{F,n}(G(t))\right) =s_{F^{\varphi},n}(G^\varphi(t))
 =\alpha(G^\varphi(t_1)\cdots G^\varphi(t_n)) =G^\varphi(\alpha_1)
 \cdots G^\varphi(\alpha_n).$$
  which is the formula of the first part of the Lemma.

The second part of the Lemma follows from the first because
$g^\varphi(x) =g(x)$ for all polynomials $g(x)$ in $k[x]$.
\enddemo

\subhead 2. Roots of $F^\varphi(x)$ and invertible elements in
$A[x]/(F(x))$\endsubhead 

\definition{2.1. Notation}We shall use the notation of Sections
(1.1)\marg and (1.2). Given a ring $A$ and a prime ideal $P$ we write
$\kappa(P) =A_P/PA_P$ for the residue field. Let $k$ be a field and
assume that $A$ is a $k$--algebra. Denote by $\local$ the localization
of $k[x]$ in the multiplicatively closed subset $k[x]\setminus (x)$ of
polynomials $g(x)$ in $k[x]$ such that $g(0)\neq 0$. We have that
$A\otimes_k \local$ is the localization of the $k[x]$--algebra
$A\otimes_k k[x]$ in the multiplicatively closed set $k[x]\setminus
(x)$.
\enddefinition

\proclaim{2.2 Lemma}Given a field extension $K$ of $k$. Let $G(x)$ be
a polynomial in $K[x]$ that has a non--zero root in the algebraic
closure $\overline{K}$ of $K$ which is algebraic over $k$. Then there
is a polynomial $g(x)$ in $k[x]$ with $g(0)\neq 0$ and a factorization
$I(x)g(x) =H(x)G(x)$ in $K[x]$, where $I(x)$ is a non--zero
polynomial with $\deg(I) < \deg(G)$ whose roots in $\overline{K}$
are zero or transcendental over $k$.
\endproclaim

\demo{Proof}We shall prove the Lemma by induction on the degree $m$ of
$G(x)$.

Let $\alpha$ be a non--zero root of $G(x)$ in $\overline{K}$ which is
algebraic over $k$. Denote by $g_1(x)\in k[x]$ and $G_1(x) \in K[x]$
the minimal polynomials of $\alpha$ over $k$ respective $K$. Then we
have that $g_1(0)\neq 0$ and $\deg(G_1) \geq 1$. Moreover we have
factorizations $g_1(x) =H_1(x) G_1(x)$ and $G(x) =I_1(x) G_1(x)$ in
$K[x]$, where $I_1(x)$ is non--zero and $\deg(I_1)
<\deg(G)$. Consequently $I_1(x)g_1(x) =H_1(x)I_1(x) G_1(x) =H_1(x)
G(x)$ in $K[x]$.

When $m=1$ we have that $I_1(x)$ is a non--zero constant in $K$ and
the Lemma holds. If all the roots of $I_1(x)$ are zero or
transcendental over $k$ we have proved the Lemma.

 Assume that the Lemma holds for
all polynomials of degree less that $m$. It remains to prove the Lemma
when $m>1$ and $I_1(x)$ has a non--zero root in $\overline{K}$ that is
algebraic over $k$. Since $\deg(I_1) =\deg(G) -\deg(G_1) <m$ it
follows from the induction assumption that there is a polynomial
$g_2(x) \in k[x]$ such that $g_2(0) \neq 0$, and a factorization $I(x)
g_2(x) =H_2(x) I_1(x)$ in $K[x]$, where $I(x)$ is a non--zero
polynomial such that $\deg(I) < \deg(I_1)$ whose roots are all zero or
transcendental over $k$. We get that $I(x) g_1(x) g_2(x) =H_2(x)
g_1(x) I_1(x) =H_1(x) H_2(x) G(x)$, and we have proved the Lemma.
\enddemo

\proclaim{2.3. Theorem}Given a $k$--algebra $A$ and a polynomial
 $$F(x) =x^n -u_1 x^{n-1} +\cdots + (-1)^n u_n$$
 in $A[x]$. The following six assertions are equivalent:
\roster
\item For all maximal ideals $P$ of $A$ with residue map
$\varphi\colon A\to \kappa(P)$, the roots of $F^\varphi(x)
=x^n-\varphi(u_1)x^{n-1} +\cdots 
  +(-1)^n \varphi(u_n)$ in the algebraic closure
 of $\kappa(P)$
  are zero or transcendental over $k$.
\item For every polynomial $g(x)$ in $k[x]$ with $g(0) \neq 0$ we have
  that $s_{F,n}(g(t))$ is invertible in $A$.
\item For every polynomial $g(x)$ in $k[x]$ with $g(0) \neq 0$ we have
  that the class of $g(x)$ in $A[x]/(F(x))$ is invertible.
\item The canonical fraction map
 $$A[x]/(F(x)) \to A \otimes_k \local/(F(x))$$
 is an isomorphism.
\item The $A$--module $A\otimes_k\local /(F(x))$ is free of rank $n$
  with a basis consisting of the classes of $1,x, \dots , x^{n-1}$.
\item For all maximal ideals $P$ of $A$ with residue map
$\varphi\colon A\to \kappa(P)$, the $\kappa(P)$--vector\-space
$\kappa(P)\otimes_k\local/(F^\varphi(x))$  is 
$n$--dimensional with a basis consisting of the classes of $1,x, \dots
, x^{n-1}$.
\endroster 
\endproclaim

\demo{Proof}Assume that assertion (1) holds. Let $\varphi\colon A\to
K$ be the residue map associated to a maximal ideal of $A$.  It
follows from Lemma (1.4)\marg that we have $\varphi(s_{F,n}(g(t)))
=g(\alpha_1)\cdots g(\alpha_n)$ in ${K}$, where $\alpha_1, \dots ,
\alpha_n$ are the roots of $F^\varphi(x)$ in the algebraic closure of
${K}$. Since $g(0)\neq 0$ we have that $g(\alpha_i)\neq 0$ both when
$\alpha_i$ is zero and when $\alpha_i$ is transcendental over
$k$. Hence $\varphi(s_{F,n}(g(t))) \neq 0$. Since this holds for all
maximal ideals of $A$ we have that $s_{F,n}(g(t))$ is not contained in
any maximal ideal of $A$ and thus is invertible in $A$. Consequently
assertion (2) holds.

It follows from Lemma (1.3)\marg  that assertion (2) implies assertion
(3). 

The canonical map of assertion (4) is the fraction map of the
$k[x]$--algebra $A[x]$ by the multiplicatively closed subset
$k[x]\setminus (x)$. Thus the fraction map is an isomorphism if and
only if the classes of the elements of $k[x]\setminus(x)$ in
$A[x]/(F(x))$ are invertible. Consequently assertions (3) and (4) are
equivalent.

Since the $A$--module $A[x]/(F(x))$ is free of rank $n$ with a basis
consisting of the classes $1, x, \dots , x^{n-1}$, we have that
assertions (4) and (5) are equivalent.

It is evident that assertion (5) implies assertion (6).

We shall prove that assertion (6) implies assertion (1). Assume that
assertion (1) does not hold. Then there exists a maximal ideal in $A$
with residue map $\varphi\colon A\to K$ such that $F^\varphi(x)$ has a
non--zero root in the algebraic closure of $K$ that is algebraic over
$k$. It follows from Lemma (2.2)\marg that there is a polynomial $g(x)
\in k[x]$ with $g(0)\neq 0$, and a factorization $I(x)g(x) =H(x)
F^\varphi(x)$ in $K[x]$ where $I(x)$ is a non--zero polynomial with
$\deg(I) < \deg(F^\varphi)=n$, such that the roots of $I(x)$ are zero
or transcendental over $k$. Since $g(x)$ is invertible in $K\otimes_k
\local$ we obtain that $(I(x)) \subseteq (F^\varphi(x))$ in
$K\otimes_k \local$ and thus a surjection
 $$K\otimes_k \local/(I(x)) \to K\otimes_k
\local/(F^\varphi(x)).\tag{2.3.1}$$ We already proved that assertion
(1) implies assertion (6). Hence we conclude that $K\otimes_k
\local/(I(x))$ is a vectorspace of dimension $\deg(I)$. Since we have
the surjection (2.3.1)\marg the dimension of the $K$--vectorspace
$K\otimes_k \local/(F^\varphi(x))$ is at most equal to $\deg(I)$ and
thus strictly smaller than $n$. Hence assertion (6) does not hold. We
have thus proved that assertion (6) implies assertion (1).

\enddemo

\proclaim{2.4. Corollary}Assume that $A=K$ is a field. Given a
polynomial $G(x)$ in $K[x]$ of degree $n$. Then the $K$--vectorspace
$K\otimes_k \local/(G(x))$ is generated by the classes of $1, x, \dots
, x^{n-1}$.
\endproclaim

\demo{Proof}If the roots of $G(x)$ in the algebraic closure
$\overline{K}$ of $K$ are zero or transcendental over $k$, the
Corollary follows from the Theorem.

Assume that $G(x)$ has a non--zero root in $\overline{K}$ that is
algebraic over $k$. It follows from Lemma (2.2)\marg that there is a
polynomial $g(x)$ in $k[x]$ with $g(0)\neq 0$, and a factorization
$I(x) g(x) =H(x)G(x)$ in $K[x]$ where $I(x)$ is a non--zero polynomial
with $\deg (I) <\deg(G)$ whose roots in $\overline{K}$ are zero or
transcendental over $k$. We obtain that $(I(x))\subseteq (G(x))$ in
$K\otimes_k \local$ and thus a surjection $K\otimes_k \local/(I(x))
\to K\otimes_k \local/(G(x))$. By Theorem (2.3)\marg we have that
$K\otimes_k \local/(I(x))$ is generated by the classes of $1, x, \dots
, x^{m-1}$ with $m=\deg(I) < \deg(G) =n$. Thus $K\otimes_k
\local/(G(x))$ is generated by the classes of $1, x, \dots , x^{n-1}$.
\enddemo

\subhead 3. Ideals in $A\otimes_k\local$ whose residue rings are free
$A$--modules\endsubhead

\definition{3.1. Notation}Given a field $k$ and a $k$--algebra $A$.
We denote by $k[x]$ the polynomial ring in a variable $x$ over $k$ and
write $A[x] =A\otimes_k k[x]$.  We denote by $\local$ the localization
of $k[x]$ in the multiplicatively closed set $k[x]\setminus (x)$
consisting of polynomials $g(x)$ of $k[x]$ with $g(0) \neq 0$.
\enddefinition

\proclaim{3.2. Lemma}Given an ideal $I$ of $A\otimes_k \local$ such
that the residue ring $A\otimes_k\local/I$ is a free $A$--module of
rank $n$. Then the $A$--module $A\otimes_k \local /I$ has a basis
consisting of the classes of the elements $1, x, \dots , x^{n-1}$.
\endproclaim

\demo{Proof}We must prove that the $A$--module homomorphism
 $$A^n \to A\otimes_k \local /I\tag{3.2.1}$$
 which sends the coordinates of $A^n$ to the classes of $1, x, \dots ,
 x^{n-1}$ is an isomorphism. It suffices to prove that the
 localization of the map (3.2.1)\marg in each prime ideal of $A$ is an
 isomorphism. Hence we may assume that $A$ is a local $k$--algebra. In
 fact, it suffices to prove that the map (3.2.1)\marg is surjective
 since any set of $n$ generators of a free module of rank $n$ form a
 basis.

Assume that $A$ is local and let $K$ be the residue field of $A$. We
denote by $I_K$ the image of the ideal $I$ by the residue map
$A\otimes_k \local \to K\otimes_k \local$. From
the map (3.2.1)\marg we obtain a homomorphism 
$$K^n \to K\otimes_k \local /I_K\tag{3.2.2}$$ 
 of $K$--vectorspaces which
sends the coordinates of $K^n$ to the classes of $1, x, \dots ,
x^{n-1}$. If the map (3.2.2)\marg is injective then it is surjective
because $\dim_K(K\otimes_k \local/I_K) =n$ by assumption. On the other
hand, if (3.2.2)\marg is not injective, there is a polynomial $G(x)
=x^m+g_1x^{m-1} +\cdots +g_m$ in $K[x]$ of degree $m\leq n$ such that
$(G(x)) \subseteq I_K$ in $K\otimes_k \local$. Hence we have a
surjection $K\otimes_k \local/(G(x)) \to K\otimes_k \local/I_K$.  It
follows from Corollary (2.4)\marg that $K\otimes_k \local/(G(x))$ is
generated by the classes of $1, x, \dots , x^{m-1}$. Thus the
$K$--vectorspace $K\otimes_k \local/I_K$ is generated by the classes of
$1, x, \dots , x^{n-1}$, and the map (3.2.2)\marg is surjective even
when it is not injective.

By assumption the $A$--module $A\otimes_k \local /I$\/ is finitely
generated. Hence it follows from Nakayama's Lemma that the $A$--module
homomorphism (3.2.1) is surjective.  \enddemo

\proclaim{3.3. Theorem}Given an ideal $I$ in $A\otimes_k \local$ such
that $A\otimes_k\local /I$ is a free $A$--module of rank $n$. Then the
ideal $I$ is generated by a unique monic polynomial
 $$F(x) =x^n -u_1 x^{n-1} +\cdots +(-1)^n u_n$$
 in $A[x]$.
\endproclaim

\demo{Proof}Since the $A$--module $A\otimes_k \local/I$\/ is free of
rank $n$ by assumption, it follows from Lemma (3.2)\marg that the
$A$--module $A\otimes_k \local/I$\/ has a basis consisting of the
classes of the elements $1, x, \dots , x^{n-1}$. Hence the class of
$x^n$ in $A\otimes_k \local /I$\/ can be written as a unique
$A$--linear combination of the classes of $1, x, \dots , x^{n-1}$. It
follows that there is a unique monic polynomial $F(x) =x^n -u_1
x^{n-1} +\cdots +(-1)^n u_n$ in $A[x]$ whose image is contained in $I$.

To show that $F(x)$ generates $I$\/ we must prove that the surjective
residue map
 $$A\otimes\local/(F(x)) \to A\otimes_k \local /I\tag{3.3.1}$$
  is injective. It suffices to prove that the localization of
(3.3.1)\marg at every prime ideal of $A$ is injective. Hence we may
assume that $A$ is local. Denote by $\varphi\colon A\to K$ the residue
map, and let $I_K$ be the image of the ideal $I$ by the residue map
$A\otimes_k \local \to K\otimes_k \local$. From the map
(3.3.1)\marg we obtain the $K$--linear residue map
 $$K\otimes_k \local/(F^\varphi(x)) \to K\otimes_k
\local/I_K.\tag{3.3.2}$$  
By assumption we have that the $A$--module $A\otimes_k \local/I$\/ is
free of rank $n$. Consequently we have that $K\otimes_k \local /I_K$\/
is an $n$--dimensional $K$--vectorspace. It follows from Lemma
(3.2)\marg that the classes of $1,x, \dots , x^{n-1}$ in $K\otimes_k
\local /I_K$\/ form a $K$--basis. From Corollary (2.4)\marg it follows
that the classes of $1, x, \dots , x^{n-1}$ in the $K$--vectorspace
$K\otimes_k\local/(F^\varphi(x))$ are generators. The existence of the
sur\-jection (3.3.2) therefore shows that the classes of $1, x, \dots
, x^{n-1}$ in $K\otimes_k \local /(F^\varphi(x))$ form a
$K$--basis. Hence it follows from Theorem (2.3)\marg that the roots of
$F^\varphi(x)$ in the algebraic closure of $K$ are zero or
transcendental over $k$. Consequently it follows from Theorem
(2.3)\marg that the classes of $1, x, \dots , x^{n-1}$ in the
$A$--module $A\otimes_k \local /(F(x))$ form a basis. On the other
hand it follows from Lemma (3.2)\marg that the classes of $1, x, \dots
, x^{n-1}$ in the $A$--module $A\otimes_k \local /I$\/ form a
basis. It follows that the map (3.3.1) is injective. We have proved
the Theorem.
\enddemo

\subhead 4. The coordinate ring of the Hilbert scheme\endsubhead

\definition{4.1. Notation}Given a field $k$. Denote by $k[x,t] = k[x,
t_1, \dots , t_n]$ the polynomial ring over $k$ in the variables $x,
t_1, \dots , t_n$. For every $k$--algebra $A$ we write $A[x,t]
=A\otimes_k k[x,t]$. We denote by $s_i(t)$ the $i$'th elementary
symmetric polynomial in the variables $t_1, \dots , t_n$. For every
polynomial $g(x)$ in $k[x]$ we form, as in the construction of Section
(1.2),\marg the symmetric polynomial $s_n(g(t)) =g(t_1)\cdots g(t_n)$
in $k[t_1, \dots , t_n]$. The set
 $$U=\{s_n(g(t)) \colon g(x) \in k[x] \text{ and } g(0) \neq 0\}$$
form a multiplicatively closed subset of $k[s_1(t), \dots ,
s_n(t)]$. We write
 $$H_n =U^{-1} k[s_1(t), \dots , s_n(t)]$$
 and 
$$F_n(x) =x^n -s_1(t) x^{n-1} +\cdots +(-1)^n s_n(t)$$ 
 in $H_n[x]$.
\enddefinition

\proclaim{4.2 Proposition}Given a $k$--algebra homomorphism
 $$\psi \colon H_n =U^{-1} k[s_1(t), \dots , s_n(t)]\to A$$
 Let $\psi(s_i(t)) =u_i$ for $i=1, \dots , n$ and write
 $$F_n^\psi(x) =x^n -u_1x^{n-1} +\cdots +(-1)^n u_n$$
in $A[x]$. Then the $A$--module 
 $$A\otimes_k\local/(F_n^\psi(x))$$
 is free of rank $n$ with a basis consisting of the classes of the
elements $1, x, \dots , x^{n-1}$.

In particular the $H_n$ module $H_n\otimes_k \local/(F_n(x))$ is
free of rank $n$ with a basis consisting of the classes of the
elements $1, x, \dots, x^{n-1}$.
\endproclaim

\demo{Proof} The fraction map $\xi\colon k[s_1(t), \dots , s_n(t)] \to
H_n$ composed with $\psi$ defines a map $\zeta\colon k[s_1(t), \dots ,
s_n(t)] \to A$. We have that $\zeta$ is the restriction to
$k[s_1(t), \dots , s_n(t)]$ of the map $u$ of the construction of
Section (1.2) for the polynomial $F_n^\psi(x)$.\marg Hence we have
that
 $$\psi (s_n(g(t))) =\zeta\left(s_n(g(t))\right) =s_{F^\psi,n}(g(t))$$
for all polynomials $g(x)$ of $k[x]$. Given a polynomial $g(x)$ of
$k[x]$ with $g(0)\neq 0$.  Since $\zeta$ factors via the fraction map
$\xi$ and $s_n(g(t))$ is invertible in $H_n$ by definition, we have
that $\zeta(s_n(g(t))) =s_{F^\psi,n}(g(t))$ is invertible in
$A$. Consequently it follows from Theorem (2.3)\marg that
$A\otimes_k\local/(F_n^\psi(x))$ is a free $A$--module with a basis
consisting of the classes of the elements $1, x, \dots , x^{n-1}$. We
have proved the Proposition.
\enddemo

\definition{4.3. Notation}We denote by $\Cal Hilb_a^n$ the {\it affine
Hilbert functor} from the category of $k$--algebras to sets that sends
a $k$--algebra $A$ to
 $$\align
\Cal Hilb_a^n(A) =& \{\text{Ideals }  I \text{ in } A\otimes_k \local
\text{ such that}\\
&A\otimes_k \local/I \text{ is a free } A
\text{--module of rank } n\}.
\endalign$$
 It follows from Proposition (4.2)\marg that we, for every
 $k$--algebra $A$, obtain a natural map
 $$\Cal F(A) \colon \Hom_k(H_n, A) \to \Cal Hilb_a^n(A)$$
 which sends a $k$--algebra homomorphism $\psi\colon H_n \to A$ to
the ideal $(F_n^\psi(x))$ in $A\otimes_k \local$.
 Clearly $\Cal F$ defines a morphism of functors from $k$--algebras
to sets. 
\enddefinition

\proclaim{4.4. Theorem}The morphism $\Cal F$ is an isomorphism of
functors.  Equivalently, the $k$--algebra $H_n$ represents the functor
$\Cal Hilb_a^n$.
\endproclaim 

\demo{Proof}We shall construct an inverse to $\Cal F$. Given a
$k$--algebra $A$ and an ideal $I$ in $A\otimes_k \local$ such that the
$A$--module $A\otimes_k\local/I$ is free of rank $n$. It follows from
Theorem (3.3)\marg that there is a unique polynomial $F(x)=x^n
-u_1x^{n-1} +\cdots +(-1)^n u_n$ in $A[x]$ such that $(F(x)) =I$ in
$A\otimes_k \local$. In particular we have that $A\otimes_k
\local/(F(x))$ is a free $A$--module of rank $n$. Hence it follows
from Theorem (2.3)\marg that the element $s_{F,n}(g(t))$ is invertible
in $A$ for all polynomials $g(x)$ in $k[x]$ such that $g(0) \neq 0$.

We have that the $k$--algebra homomorphism $\zeta\colon k[s_1(t),
\dots , s_n(t)]\to A$ determined by $\zeta(s_i(t))=u_i$ for $i=1,
\dots , n$ coincides with the restriction to $k[s_1(t), \dots ,
s_n(t)]$ of the homomorphism $u$ of Section (1.2).\marg It follows
that $\zeta(s_i(g(t))) =s_{F,i}(g(t))$ for $i=1, \dots , n$ and for
all polynomials $g(x)$ in $k[x]$. In particular the elements
$s_n(g(t))$ of $U$ are mapped to the invertible elements
$s_{F,n}(g(t))$ in $A$. Consequently the $k$--algebra homomorphism
$\zeta$ factors through a unique $k$--algebra homomorphism $\psi\colon
H_n = U^{-1} k[s_1(t), \dots , s_n(t) ] \to A$. Thus we have
constructed a map $\Cal Hilb_a^n(A) \to \Hom_k(H_n, A)$. It is easy to
check that this map is the inverse to $\Cal F(A)$.
\enddemo

\definition{4.5. Definition}We define the {\it Hilbert functor $\Cal
Hilb^n$ of families of length $n$ subschemes of $\Spec \local$} as the
contravariant functor from schemes over $k$ to sets, that to a
$k$--scheme $T$ associates the set
 $$\align 
\Cal Hilb^n(T) =&\{\text{Closed subschemes } Z \text{ in }
T\times_{\Spec k} \Spec \local \text{ such that }\\ &p_\ast\Cal O_Z
\text{ is a locally free } \Cal O_T \text{--module of rank } n \text{,
where}\\ &p\colon T\times_{\Spec k} \Spec \local \to T \text{ is the
projection}\}.
\endalign$$
 Equivalently the set $\Cal Hilb^n(T)$ consists of the quasi--coherent
 ideals $\Cal I$ in $\Cal O_T\otimes_{\Cal O_{\Spec k}}\Cal O_{\Spec
 \local}$  such that $\left(\Cal O_T\otimes_{\Cal O_{\Spec k}} \Cal
 O_{\Spec \local}\right)/\Cal I$ is locally free considered as an $\Cal
 O_T$--module via $p$.
\enddefinition

\proclaim{4.6. Theorem}We have that the scheme $\Spec H_n$ with
 the universal family $\Spec \left(H_n\otimes_k \local /(F_n(x))\right)$
 represents $\Cal Hilb^n$.
\endproclaim

\demo{Proof} For every $k$--scheme $T$ and every point $Z$ in $\Cal
 Hilb^n(T)$ we cover $T$ with affine open subsets $\Spec A$ such that
 $p_\ast\Cal O_Z\vert \Spec A$ is a free $\Cal O_{\Spec
 A}$--module. It follows from Theorem (4.4)\marg that we have a unique
 morphism $\Spec A\to \Spec H_n$ such that the family $\Spec \left(
 H_n\otimes_k \local /(F_n(x))\right)$ pulls back to the family $Z\cap
 p^{-1}(\Spec A)$ over $\Spec A$. Since the map is unique the
 morphisms for the affine subsets covering $T$\/ glue together to a
 morphism $T\to \Spec H_n$ such that the family $\Spec \left( H_n
 \otimes_k \local/(F_n(x))\right)$ pulls back to $Z$.
\enddemo

\definition{4.7. Note}The point $s_1(t)=\cdots =s_n(t)=0$ is, as
expected, the only $k$--rational point of $\Spec H_n$. Indeed, let
$(u_1, \dots , u_n)$ be point different from the origin of
the affine $n$--dimensional space $\bold A_k^n$ over $k$.  To show
that this point does not lie in $\Spec H_n$ we chose roots $\alpha_1,
\dots , \alpha_n$ in $\overline{k}$ of the polynomial $x^n
-u_1 x^{n-1} +\cdots +(-1)^n u_n$ in $k[x]$. Since all
the $u_i$ are not zero there is at least one root $\alpha_i$ which
is not zero. Let $g(x) \in k[x]$ be the minimal polynomial of such a
root. Then $g(0) \neq 0$ and thus $s_n(g(t))\in U$. The map
$\alpha\colon k[t_1, \dots , t_n] \to \overline k$ which sends $t_i$ to
$\alpha_i$ for $i=1, \dots , n$ iduces the map $u\colon k[s_1(t),
\dots , s_n(t)] \to k$ which sends $s_i(t)$ to $u_i$ for $i=1,
\dots ,n$. However, the map $u$ does not factor through the
fraction map $k[s_1(t), \dots , s_n(t)] \to H_n$ because $s_n(g(t))
=g_1(t)\cdots g_n(t)$ is mapped to the element $g(\alpha_1)\cdots
g(\alpha_n) =0$ in $k$. Hence the point $(u_1, \dots , u_n)$
is not in $\Spec H_n$.

When $n=1$ we have that $H_1=\local$ is a local ring. However,
when $n\geq 2$ the ring $H_n$ is not local. To see this we first
prove that the ideal $(F(t))$ in $k[s_1(t), \dots , s_n(t)]$ generated
by a non--constant, symmetric polynomial $F(t)$ which is irreducible
in $k[t_1, \dots , t_n]$ does not intersect the multiplicatively
closed subset $U$ of $k[s_1(t), \dots , s_n(t)]$. Assume that $(F(t))$
intersects $U$. Then we have that $F(t)G(t) =s_n(f(t)) =f(t_1) \cdots
f(t_n)$ in $k[s_1(t), \dots , s_n(t)]$ for a polynomial $f(x) \in
k[x]$ with $f(0)\neq 0$. Then $F(t)$ divides one of the polynomials
$f(t_1), \dots , f(t_n)$ in $k[t_1, \dots , t_n]$. Hence $p(t)$ is a
polynomial in one of the variables $t_i$. Hence, when $n\geq 2$ it can
not be symmetric, contrary to our assumption. Thus $(F(t))$ does not
intersect $U$.

For each non--constant, symmetric, irreducible polynomial $F(t)$ in
$k[t_1, \dots , t_n]$ we can choose an ideal $P_F$ which is maximal
among the ideals in $k[s_1(t), \dots , s_n(t)]$ that contain $F(t)$
and do not intersect $U$. Then $P_F$ is a prime ideal and the ideal
$P_FH_n$ is maximal in $H_n$. In this way we can construct an
abundance of maximal ideals in $H_n$ when $n\geq 2$. For example we
can choose $F_v(t) =v+s_1(t)$ with $v\in k$. It is clear that the
maximal ideals $P_{F_v}H_n$ for different $v$ are all different.

\enddefinition


\Refs
\tenpoint
\widestnumber\key{ABC}

\ref \key B
\by J. Brian\c{c}on
\paper {Description de $\operatorname{Hilb}^n{\bold C} \{x,y\}$}
\jour {Invent. Math.}
\pages {45--89}
\vol {41}
\issue {1}
\yr {1977}
\endref

\ref \key BI
\by {J. Brian\c{c}on and A. A. Iarrobino}
\paper {Dimension of the punctual {H}ilbert scheme}
\jour {J. Algebra}
\vol {55}
\issue {2}
\pages {536--544}
\yr {1978}
\endref

\ref \key C
\by {M. Coppens}
\paper {The fat locus of {H}ilbert schemes of points}
\jour {Proc. Amer. Math. Soc.}
\pages {777--783}
\vol {118}
\issue {3}
\yr {1993}
\endref

\ref \key G
\by {M. Granger}
\paper {G\'{e}om\'{e}trie des sch\'{e}mas de {H}ilbert ponctuels}
\jour {M\'{e}m. Soc. Math. France. (N.S.)}
\issue {8}
\pages {84}
\yr {1983}
\endref

\ref \key I1
\by {A. A. Iarrobino}
\paper {Punctual {H}ilbert {S}chemes}
\jour {Mem. Amer. Math. Soc.}
\pages {viii+112}
\vol {10} 
\issue {188}
\yr {1977}
\endref

\ref \key I2
\by {A. A. Iarrobino}
\paper {Punctual {H}ilbert {S}chemes}
\jour {Bull. Amer. Math. Soc.}
\vol {78}
\pages {819--823}
\yr {1972}
\endref

\ref \key I3
\by {A. A. Iarrobino}
\paper {Hilbert {s}cheme of {p}oints: overview of last ten years}
\inbook {Algebraic {G}eometry, Bowdoin 1985, Proc.  Sympos. Pure Math.}
\vol {46, part 2}
\pages {297--320}
\yr {1987}
\publ {Amer. Math. Soc}
\publaddr {Providence, RI}
\endref

\ref \key P
\by {G. Paxia}
\paper {On flat families of fat points}
\jour {Proc. Amer. Math. Soc.}
\pages {19--23}
\vol {112 }
\issue {1}
\yr {1991}
\endref

\ref \key S1
\by {R.M. Skjelnes}
\paper {On the representability of $\Cal Hilb^n k[x]_{(x)}$}
\yr {1999}
\paperinfo {To appear in J. London Math. Soc.}
\endref

\ref \key S2
\by {R.M. Skjelnes}
\paper {Symmetric tensors with applications to Hilbert schemes}
\yr {1999}
\paperinfo {To appear}
\endref


\endRefs
\enddocument

\end